\documentstyle[12pt]{article}

\begin{document}

\title{Canonical Generators of the Cohomology of Moduli
of Parabolic Bundles on Curves}
\author{I.Biswas ~~and~~ N.Raghavendra\thanks{On leave
from the School of Mathematics, SPIC Science Foundation,
Madras, India.}\\
{\small Mathematics Section, International Centre for
Theoretical Physics}\\ {\small P.O.Box 586, 34100 Trieste, Italy.}\\
{\small E-mail:~~biswas@ictp.trieste.it ~and~raghu@ictp.trieste.it}\\
{\small Fax: ~39 40 224163}}
\date{}
\maketitle

\newtheorem{thm}{Theorem}[section]
\newtheorem{cor}[thm]{Corollary}
\newtheorem{lemma}[thm]{Lemma}
\newtheorem{propn}[thm]{Proposition}
\newtheorem{rem}[thm]{Remark}
\newtheorem{assume}[thm]{Assumption}
\newtheorem{nota}[thm]{Notation}
\newtheorem{defn}[thm]{Definition}
\newtheorem{exam}[thm]{Examples}

\newcommand{\ay}{\mbox{${\cal A}$}}
\newcommand{\uxnd}{\mbox{${\cal U}_{X}(n,d)$}}
\newcommand{\lra}{\mbox{$\longrightarrow$}}
\newcommand{\tsux}{\mbox{${\cal SU}_{X}(2,L,\Delta)$}}
\newcommand{\e}{\mbox{${\cal E}$}}
\newcommand{\ite}{\mbox{$i_{t}$}}
\newcommand{\jay}{\mbox{$j_{x}$}}
\newcommand{\p}{\mbox{${\bf P}$}}
\newcommand{\q}{\mbox{${\bf Q}$}}
\newcommand{\g}{\mbox{${\cal G}$}}
\newcommand{\asp}{\mbox{${\cal A}^{s}_{\rm par}$}}
\newcommand{\be}{\mbox{${\cal B}$}}
\newcommand{\ep}{\mbox{${\cal E}_{\rm par}$}}
\newcommand{\gp}{\mbox{${\cal G}_{\rm par}$}}
\newcommand{\gbp}{\mbox{$\bar{\cal G}_{\rm par}$}}
\newcommand{\ebp}{\mbox{$\bar{\cal E}_{\rm par}$}}
\newcommand{\eb}{\mbox{$\bar{\cal E}$}}
\newcommand{\gb}{\mbox{$\bar{\cal G}$}}
\newcommand{\ux}{\mbox{${\cal U}_{X}(n,d,\Delta )$}}
\newcommand{\sux}{\mbox{${\cal S}{\cal U}_{X}(n,L,\Delta )$}}
\newcommand{\jx}{\mbox{$J_{X}$}}
\newcommand{\gn}{\mbox{$\Gamma _{X}(n)$}}
\newcommand{\homo}{\mbox{${\cal H}{\it om}$}}
\newcommand{\exi}[1]{\mbox{${\cal E}_{\rm par}^{x,#1}$}}
\newcommand{\uxi}[1]{\mbox{$U^{x,#1}$}}
\newcommand{\en}{\mbox{${\cal E}{\it nd}$}}
\newcommand{\ra}{\mbox{$\rightarrow$}}
\newcommand{\jxs}{\mbox{$j_{x}^{\ast}$}}
\newcommand{\its}{\mbox{$i_{t}^{\ast}$}}
\newcommand{\f}{\mbox{${\cal F}$}}


\section{Introduction}

The aim of this paper is to determine generators of the
rational cohomology algebras of moduli spaces of parabolic
vector bundles on a curve, under some `primality' conditions
(see Assumptions 1.1 and 1.2)
on the parabolic datum. These generators are
canonical in a sense which will be made precise below.
Our results are new even for usual vector bundles (i.e., vector
bundles without parabolic structure) whose rank is greater
than $2$
and is coprime to the degree; in this case, they are
generalizations of a theorem of Newstead \cite{N}, where the
case of vector bundles of rank $2$ and odd degree is studied.

Let $X$ be a compact Riemann surface, fix integers $n$ and $d$
with $n$ positive, and let $\Delta $ be a parabolic datum of
rank $n$ on $X$ (see Section 3 below).
Denote by $\ux $ the moduli space of parabolic vector bundles
of rank $n$ and degree $d$, which are parabolic semistable with
respect to $\Delta $. Fix a holomorphic line bundle $L$ of
degree $d$ on $X$, and let $\sux $ be the subvariety
of $\ux $ consisting of vector bundles with determinant
isomorphic to $L$. We make the following
two hypotheses on the parameters $n$, $d$ and $\Delta$.

\begin{assume}
\begin{em}
Every parabolic vector bundle of rank $n$ and degree $d$ on
$X$ which is parabolic semistable with respect to the
parabolic datum $\Delta $ is in fact parabolic stable.
\end{em}
\end{assume}

\begin{assume}
\begin{em}
There exists a universal parabolic bundle (or briefly, a
universal bundle) on $\ux \times X$.
\end{em}
\end{assume}

Recall that a universal bundle on $\ux \times X$ is a vector
bundle $U$ on $\ux \times X$ together with a flag of
subbundles $\jxs U = \uxi{1} \supset \uxi{2} \supset \cdots
\supset \uxi{k_{x}} \supset \uxi{k_{x}+1} = 0$ in $\jxs U$ for
each $x \in J$ ($J$ being the set of parabolic points), where
$\jay : \ux \ra \ux \times X$ is the map $E \mapsto (E,x)$,
such that for each $E \in \ux$, the restriction of $U$ to
$\{ E \} \times X$ is parabolically isomorphic to $E$. We use
the same symbols to denote the restrictions of $U$ and
$\uxi{i}$ etc. to $\sux $.

\begin{nota}

\begin{em}
\begin{itemize}
\item[]
\item If $S$ and $T$ are topological spaces, and if $\alpha
\in H^{\ast }(S\times T,\q)$, then $\sigma (\alpha
):H_{\ast}(T,\q) \rightarrow H^{\ast} (S,\q)$ is the map
$\sigma (\alpha )z=\alpha /z$, where $/$ denotes the slant
product.\\
\item If $V$ is a vector bundle, then $\p(V)$ denotes its
projectivization.\\
\item All cohomologies in the paper have  $\q$-coefficients.
\end{itemize}
\end{em}
\end{nota} \label{1.3}

Having settled on the notation, we now have the following two
theorems which are the main results of this paper.

\begin{thm}
Suppose that Assumptions 1.1 and 1.2 hold, and let $U$ be any
universal bundle on $\ux \times X$. Then, the rational
cohomology algebra of $\ux $ is generated by the Chern classes
$c_{j}(\homo (U^{x,i} ,U^{x,i-1} ))~(x\in J)$ and the images of
$$\sigma (c_{1}(U)): H_{1}(X) \rightarrow H^{1}(\ux )~~~{\it
and}$$
$$\sigma (a_{i}(\p (U))): H_{r}(X) \rightarrow H^{2i-r}(\ux
)~~(2\leq i \leq n,~0\leq r \leq 2)$$
where $a_{i}(.)$ denote the characteristic classes of
projective bundles introduced in Definition 2.4 below.
\end{thm}

\begin{thm}
Suppose that Assumptions 1.1 and 1.2 hold, and let $U$ be any
universal bundle on $\sux \times X$. Then the rational
cohomology algebra of $\sux $ is generated by the Chern
classes $c_{j}(\homo (U^{x,i},U^{x,i-1}))~(x\in J)$
and the images of
$$\sigma (a_{i}(\p (U))): H_{r}(X) \rightarrow H^{2i-r}(\sux
)~~(2\leq i\leq n,~0\leq r \leq 2).$$
\end{thm}

Note that the generators given in Theorems 1.4 and
1.5 are
{\em canonical}, i.e., independent of the choice of a
universal bundle (which is easily seen to be non-unique).
Indeed, if $U'$ is another universal bundle, then there exists
a line bundle $\xi$ on $\ux$ such that
$U' \cong U\otimes p^{*}\xi$,
where $p : \ux \times X \ra \ux$ is the canonical projection.
Now, for every $z\in H_{1}(X)$, we have $\sigma (c_{1}(U'))z =
\sigma (c_{1}(U))z$, since $\sigma (c_{1}(p^{*}\xi ))z = 0$;
on the other hand, it is obvious that $\p (U') \cong \p (U)$
and $\homo ((U')^{x,i},(U')^{x,i-1}) \cong
\homo (\uxi{i}, \uxi{i-1})$.

We now relate the above theorems to certain results of Atiyah
and Bott \cite{AB}. Let $\uxnd$ be the moduli space of stable
vector bundles of rank $n$ and degree $d$, where $n$ and $d$
are coprime. In this case, Atiyah and Bott \cite{AB} proved
that the Kunneth components (with respect to any basis of
$H^{*}(X,\q )$) of the Chern classes of any universal bundle
on $\uxnd \times X$ generate the rational cohomology algebra
of $\uxnd$. Theorem 1.4 above differs from this result
in the following respects. Firstly, we work throughout in the
setup of parabolic bundles, whereas Atiyah and Bott were
working with usual vector bundles. Secondly, as we have
observed above, the generators we
obtain are canonical, i.e., they are independent of the choice
of a universal bundle, whereas the Kunneth components of the
Chern classes of a universal bundle $U$ do depend on the
choice of $U$. Finally, by specializing to the case where the
parabolic set is empty, and applying Lemma 2.6 below, we
obtain the above result of Atiyah and Bott from Theorem
1.4; whereas it does
not seem possible to deduce Theorem 1.4 from the result
of Atiyah and Bott: the difficulty is due to the fact that the
slant product does not behave well with the cup product.

We should remark that Beauville \cite{Beau} has given
another proof of the above result of Atiyah and Bott. In the
parabolic setup, and under Assumptions 1.1 and 1.2, the
method of
Beauville can be used to deduce that the
Kunneth components of the
Chern classes of any universal bundle $U$ and the Chern classes
$c_{j}(\homo (\uxi{i},\uxi{i-1}))$ generate the cohomology of
$\ux$, a statement which, as we have seen above, is weaker
than Theorem 1.4.

The following is a consequence of the above theorems.

\begin{sloppy}
\begin{cor}
Let n=2,
suppose Assumptions 1.1 and 1.2 are true,
and let $U$  be
any universal bundle on $\tsux \times X$.
Then, the rational cohomology
algebra of $\tsux $ is generated by the Chern classes
$c_{j}(\homo (S^{x}, Q^{x}))$ ($x \in J$), and the image of
$$\sigma(c_{2}(\en U)) : H_{r}(X) \lra H^{4-r}(\sux )
{}~~~(0 \leq r \leq 2),$$
where $\jxs U = \uxi{1} \supset \uxi{2} = S^{x} $ is the
flag in $\jxs U$ ($x\in J$), and
$Q^{x} = \jxs U / S^{x}$.
\end{cor}
\end{sloppy}

The above corollary is a generalization to parabolic bundles
of a theorem of Newstead \cite{N}.

In our approach, Assumption 1.1 is natural and indispensable;
it is a technical necessity which guarantees that the action
of a certain gauge group on a certain space of holomorphic
structures is free. Granted this, Assumption 1.2 is not too
stringent a restriction, as the following observation shows.

\begin{propn}
Suppose the parameters $n$, $d$ and $\Delta$ satisfy
Assumption 1.1. Then, they satisfy Assumption 1.2 if any one
of the following three conditions is satisfied:
\begin{itemize}
\item The rank $n$ and the degree $d$ are coprime.
\item There exists a parabolic point $x \in J$ such that
$\sum_{i=j}^{k_{x}}n_{x,i}$
is coprime to $n$ for some $j$ ($1
\leq j \leq k_{x}$), where $n_{x,1},\ldots ,n_{x,k_{x}}$
denote the parabolic multiplicities at $x$.
\item There exists a parabolic point $x \in J$ such that
$\sum_{i=j}^{k_{x}}n_{x,i}$
and $n+d$ are coprime for some $j$.
\end{itemize}
\end{propn}

Here is a brief outline of the contents of the paper. In
Section 2, we define the characteristic classes of projective
bundles which occur in the statements of the theorems. The
next section contains a description of generators of the
rational cohomology of the classifying spaces of certain gauge
groups. The final section contains proofs of the above
results.
\section{Projective Bundles}

This preliminary section deals with some universal aspects of
projective bundles, the aim being to define explicit
characteristic classes for these bundles.

If $G$ is a topological group, then $EG \rightarrow BG$ will
denote a universal principal $G$-bundle. Cohomology groups
will have rational coefficients throughout, unless otherwise
indicated. Fix a positive integer $n$.

The natural epimorphism $\pi :U(n) \rightarrow PU(n)$ induces
a fibration $B\pi :BU(n) \rightarrow BPU(n)$ with fibre
$BU(1)$. Let $x_{1},\ldots ,x_{n}$ be the Chern roots of
$EU(n)$, so $H^{\ast}(BU(n))$ is the algebra $S[x_{1},\ldots
,x_{n}]$ of symmetric polynomials in the $x_{i}$ with rational
coefficients (or, equivalently, $H^{\ast }(BU(n))$ is the
polynomial algebra $\q [c_{1},\ldots ,c_{n}]$, where
$c_{i}=c_{i}(EU(n))$). The Leray-Hirsch theorem implies that
the map $(B\pi )^{\ast}:H^{\ast}(BPU(n))\rightarrow
H^{\ast}(BU(n))$ is injective, and
(see \cite{BH}, Section 15.2)
its image equals the subalgebra $I[x_{1},\ldots ,x_{n}]$
of $S[x_{1},\ldots ,x_{n}]$ consisting of symmetric polynomials
invariant under the affine change of variables $x_{i} \mapsto
x_{i}+d$, where $d$ is an indeterminate.

\begin{rem} \begin{em}
The above fact means that if $E\rightarrow  M$ is a vector
bundle, the characteristic classes of its projectivization $\p
(E)$ are precisely the characteristic classes of $E$ which are
invariant under tensoring by a line bundle.
\end{em} \end{rem}

\begin{sloppy}
\begin{lemma}
The above al\-gebra $I[x_{1}\ldots ,x_{n}]$ is a poly\-no\-mial
al\-gebra $\q [z_{2},\ldots ,z_{n}]$, where
$$ z_{k} = \sum _{1\leq j_{1}<\ldots ,j_{k}\leq
n}y_{j_{1}}\ldots y_{j_{k}},~~~~2\leq k \leq n,$$
are the elementary symmetric functions in $y_{1},\ldots
,y_{n}$, where
$$y_{i}=nx_{i}-\sum _{j=1}^{n}x_{j},~~~~1\leq i \leq n.$$
\end{lemma}
\end{sloppy}

The proof of the lemma is quite easy, and we omit it. Note
that the first elementary symmetric polynomial in the $y_{i}$,
namely their sum, is zero. The following assertion follows
immediately from Lemma 2.2.

\begin{cor}
For $k=2,\ldots ,n$, define $a_{k} \in H^{2k}(BPU(n))$ by
$(B\pi )^{\ast }a_{k}=z_{k}$. Then $H^{\ast }(BPU(n))$ is the
polynomial algebra $\q [a_{2},\ldots ,a_{n}]$.
\end{cor}

\begin{defn}
\begin{em}
Let $P$ be a principal $PU(n)$-bundle on a CW-complex $M$.
Then, for $k=2,\ldots ,n$, the {\em $k$-th characteristic
class} $a_{k}(P)$ of $P$ is, by definition, the element
$f^{\ast }a_{k} \in H^{2k}(M)$, where $f:M \rightarrow BPU(n)$
is some classifying map for $P$ and $a_{k}$ is as in Corollary
1.3.
(As usual, $a_{k}(P)$ is
independent of the choice of $f$.)
\end{em}
\end{defn}

\begin{exam}
\begin{em}
\begin{enumerate}
\item[] ~~~
\item[1.] If $E\rightarrow M$ is a vector bundle of rank 2,
with
Chern roots $x_{1},x_{2}$, then
$$ a_{2}(\p (E))= -(x_{1}-x_{2})^{2}=c_{2}(\en ~E),$$
where $\p(E)$ is the projectivization and $\en ~ E$ is the
endomorphism bundle of $E$.\\
\item[2.] Let $E\rightarrow M$ be a vector bundle of rank 3
such
that $c_{i}(E)=0~(i=1,2),~c_{3}\not= 0$, and let
$x_{1},x_{2},x_{3}$ be the Chern roots of $E$. Then
$$a_{3}(\p(E))=27x_{1}x_{2}x_{3} =27c_{3}(E) \not= 0,$$
while $c_{1}(\en ~E)=c_{3}(\en~E)=0$, $\en~ E$ being the
complexification of  a real vector bundle.
\end{enumerate}
\end{em}
\end{exam}

The above examples illustrate the fact that any characteristic
class of $\en~ E$ can be written as a polynomial in
$a_{i}(\p(E))$, but the converse is not true, for $\p(E)$ has
more characteristic classes than $\en ~E$.

\begin{lemma}
If $n\geq k\geq 2$, there exist a polynomial
$P_{n,k}(T_{1},\ldots ,T_{k-1})$ with rational coefficients,
and a non-zero constant $\lambda _{n,k} \in \q$ such that for
every vector bundle $E$ of rank $n$ on a CW-complex $M$, we
have
\begin{eqnarray}
a_{k}(\p(E)) = P_{n,k}(c_{1}(E),a_{2}(\p(E)),\ldots
,a_{k-1}(\p(E))) + \lambda _{n,k}c_{k}(E).
\end{eqnarray}
\end{lemma}

{\em Proof}.~~It suffices to find $P_{n,k}$ and non-zero
$\lambda _{n,k}$ such that (1) holds for $EU(n)$. Let
$u_{i}=a_{i}(\p(EU(n)))$ and $c_{i}=c_{i}(EU(n))$; then
$u_{i}=(BU)^{\ast }a_{i}$. We prove the result only for $k=2$,
the general case following easily by induction on $k$.
Since $c_{1}^{2}$ and $c_{2}$ generate $H^{4}(BU(n))$, we can
find $\alpha ,\lambda \in \q$ such that $u_{2}=\alpha
c_{1}^{2}+\alpha c_{2}$. If $\lambda =0$, then
$u_{2}-\alpha c_{1}^{2}=0$; since $\{1,c_{1},c_{1}^{2},\ldots
\}$ is an $H^{\ast }(BPU(n))$-basis of $H^{\ast}(BU(n))$, this
implies that $u_{2}=0$, contradicting the injectivity of
$(B\pi )^{\ast}$.~~$\Box$

\begin{lemma}
Let $1\leq r \leq 2$, and fix a $C^{\infty }$ vector bundle
$E$ of
rank $n$ on the $r$-sphere $S^{r}$. Then, for each $k=2,\ldots
,n$, there exists a $C^{\infty }$ vector bundle $V$ on
$S^{2k-r} \times S^{r}$ such that:
\begin{enumerate}
\item For each $t \in S^{2k -r}$, we have $V_{t} \simeq E$,
where $V_{t} = \ite ^{\ast}V,~\ite :S^{r}\ra S^{2k-r}\times
S^{r}$ being the map $x\mapsto (t,x).$\\
\item For each $x\in S^{r},~\jay ^{\ast}V$ is trivial, where
$\jay :S^{2k-r} \ra S^{2k-r}\times S^{r}$ is the map $t\mapsto
(t,x).$\\
\item $a_{k}(\p (V)) \not= 0$.
\end{enumerate}
\end{lemma}

{\em Proof}.~~ The existence of $V$ satisfying (1) and (2)
with $c_{k}(V) \not= 0$ follows by standard arguments of
K-theory. If $a_{k}(V)=0$, then Lemma 2.6 implies that
$c_{k}(V)$ is a multiple of $c_{1}(V)^{k}$, which is zero
since
$c_{1}(V)$ is
the pull-back of $c_{1}(E)$ by the second
projection, a contradiction.~~$\Box$

\section{Cohomology of Some Classifying Spaces}

In this section we describe, for later use, generators of the
rational cohomology algebras of certain mapping spaces. The
computations here are motivated by Section 2  of \cite{AB},
\cite{Dre}, and Section 5.1 of \cite{DK}. All
CW-complexes here are assumed to be finite, and cohomologies
are over $\q$.

Let $M$ be a pointed CW-complex, and let $E$ be a complex
vector bundle of rank $n$ over $M$. Fix a base point $b_{0}
\in BU(n)$. Let $\g(M,E)$ denote the complex gauge group of
$E$, with the compact-open topology. Then as shown in
\cite{AB} and \cite{Dre}, the space ${\rm
Map}_{E}(M,BU(n))$ of all maps $f:M\ra BU(n))$ such that
$f^{\ast }EU(n) \simeq E$ is a classifying space for $\g
(M,E)$, so let us denote ${\rm Map}_{E}(M,BU(n))$ by
$B\g(M,E)$. The subspace ${\rm Map}_{E}^{\ast} (M,BU(n))$
consisting of all pointed maps is closed in $B\g (M,E)$. We
denote ${\rm Map}_{E}^{\ast}(M,BU(n))$ by $\be (M,E)$. When
there is no scope for confusion, we shall write just $\g $ for
$\g (M,E)$, $\be $ for $\be (M,E)$, etc.

Let $\varepsilon :B\g \times M \ra BU(n)$ denote the evaluation
map, $(f,x) \mapsto f(x)$. Then the bundle $\varepsilon
^{\ast}EU(n)$ is called the {\it universal bundle} on $B\g
\times M$, and is denoted $\e (M,E)$. We denote the
restriction of $\e (M,E)$ to $\be \times M$ also by $\e (M,E)$.

If $\phi :M'\ra M$ is a pointed map of CW-complexes, and if
$E$ is a vector bundle over $M$, then, denoting $E'=\phi
^{\ast}E$, there is a natural map $\phi ^{\#}:B\g (M,E)\ra B\g
(M',E')$, $f\mapsto f\circ \phi$, carrying $\be (M,E)$ into $\be
(M',E')$. These constructs have the following functorial
property.

\begin{propn}
Suppose $\phi :M\ra M'$ is a pointed map of CW- complexes,
let $E$ be a complex vector bundle over $M$, and let $E'=\phi
^{\ast }E$. Write $\g $ for $\g(M,E)$, $\g '$ for $\g (M',E')$,
etc. If $d(.)$ denotes $c_{k}(.)$ or $a_{k}(\p (.))$ for some
$k$, then the diagram
\begin{center}
\setlength{\unitlength}{0.5pt}
\begin{picture}(288,210)(0,-20)

\put(20,2){\makebox(0,0){$H^{*}(B\g ')$}}
\put(256,2){\makebox(0,0){$H^{*}(B\g )$}}
\put(256,162){\makebox(0,0){$H_{*}(M)$}}
\put(20,162){\makebox(0,0){$H_{*}(M')$}}

\put(74,2){\vector(1,0){128}}
\put(256,144){\vector(0,-1){124}}
\put(74,162){\vector(1,0){128}}
\put(20,144){\vector(0,-1){124}}

\put(138,-15){\makebox(0,0){$H^{*}(\phi ^{\# })$}}
\put(304,82){\makebox(0,0){$\sigma (d(\e ))$}}
\put(138,179){\makebox(0,0){$\phi _{*}$}}
\put(-28,82){\makebox(0,0){$\sigma (d(\e '))$}}

\end{picture}
\end{center}
\noindent commutes, where $\sigma (\alpha )z=\alpha /z$
is the slant product with $\alpha $.
\end{propn}

{\em Proof.}~~The commutativity of the diagram

\begin{center}

\setlength{\unitlength}{0.5pt}
\begin{picture}(288,210)(0,-20)

\put(20,2){\makebox(0,0){$B\g \times M$}}
\put(256,2){\makebox(0,0){$BU(n)$}}
\put(256,162){\makebox(0,0){$B\g ' \times M'$}}
\put(20,162){\makebox(0,0){$B\g \times M'$}}

\put(74,2){\vector(1,0){128}}
\put(256,144){\vector(0,-1){124}}
\put(74,162){\vector(1,0){128}}
\put(20,144){\vector(0,-1){124}}

\put(138,-15){\makebox(0,0){$\varepsilon$}}
\put(286,82){\makebox(0,0){$\varepsilon '$}}
\put(138,179){\makebox(0,0){$\phi ^{\# } \times 1$}}
\put(-18,82){\makebox(0,0){$1 \times \phi$}}

\end{picture}
\end{center}
\noindent taken with the functorial properties of
the slant product, leads directly to the desired
conclusion.
{}~~$\Box $

The space $\be (M,E)$ has the following semi-universal
property.

\begin{propn}
Suppose $E$ is a complex vector bundle on a pointed
CW-complex, and $V$ a vector bundle on $T\times M$ such that:

\begin{itemize}
\item For each $t \in T$, we have $V_{t}\cong E$, where
$V_{t}=\its V$, $\ite :M\ra T\times M$ being the map $x\mapsto
(t,x)$.\\
\item If $x_{0}$ is the base point of $M$, and if
$j_{x_{0}}:T\ra T\times M$ is the map $t\mapsto (t,x_{0})$,
then $j_{x_{0}}^{*}V$ is trivial.
\end{itemize}

\noindent Then, there exists a map $\psi :T\ra \be (M,E)$ such
that
$(\psi \times 1)^{\ast}\e (M,E) \cong V$.
\end{propn}

{\em Proof}.~~General properties of classifying spaces give a
map $\theta :(T\times M, T\times \{x_{0}\})\ra (BU(n),b_{0})$
such that $\theta ^{*}EU(n)\cong V$, where $b_{0}$ is the base
point in $BU(n)$. If we define $\psi :T\ra \be $ by $\psi
(t)x=\theta (t,x)$, then $\varepsilon \circ (\psi \times
1)=\theta $, proving that $(\psi \times 1)^{*}\e \cong V$.
{}~~$\Box $

Let us study the space $\be (M,E)$ a little more when
$M=S^{r}$, $r=1,2$. Take $M=S^{1}$ first.

\begin{propn}
Let $E$ be a (necessarily trivial) complex vector bundle of
rank $n$ on $S^{1}$. Then the cohomology algebra of $\be
(S^{1},E)$ is generated by $c_{1}(\e)/[S^{1}]$ and
$a_{i}(\p(\e))/[S^{1}]$, $i=2,\ldots ,n$, where $[S^{1}]$
denotes the fundamental class of $S^{1}$.
\end{propn}

{\em Proof}.~~ Note first that $\be (S^{1},E)$ equals $\Omega
BU(n)$, which is homotopically equivalent to $U(n)$. Thus
$H^{*}(\be )$ is an exterior algebra on $n$ generators
$\theta _{i} \in H^{2i-1}(\be )$, $i=1,\ldots ,n$. Introduce
the ad-hoc notation $\omega _{1}=c_{1}(\e)/[S^{1}]$ and
$\omega _{i}=a_{i}(\p (\e ))/[S^{1}]$, $i=2,\ldots ,n$.  Since
$H^{1}(\be )=\q.\theta _{1}$, we can write $\omega
_{1}=\lambda \theta _{1}$ for some $\lambda \in \q$. Choose a
vector bundle $V$ on $S^{1}\times S^{1}$ such that
$c_{1}(V)\not= 0$. By Proposition 3.2, there exists a map
$\psi :S^{1}\ra \be $ such that $(\psi \times 1)^{*}\e = V$,
hence $c_{1}(V)/[S^{1}] = \psi ^{*}\omega _{1}= \lambda \psi
^{*}\theta _{1}$. Since $c_{1}(V)\not= 0$, this implies that
$\lambda \not=0$, and $\theta _{1}=\lambda ^{-1}\omega _{1}$.
Now let $2\leq k \leq n$ and assume that for each $i=1,\ldots
,k-1$, $\theta _{i}$ is a polynomial in $\omega _{1},\ldots
,\omega _{i}$. Thus $H^{*}(\be)$ is generated by $\omega
_{1},\ldots ,\omega _{k-1},\theta _{k},\ldots ,\theta _{n}$.
Write
\begin{eqnarray}
\omega _{k} =
P(\omega _{1},\ldots ,\omega _{k-1}) + \mu \theta _{k},
\end{eqnarray}
where $P$ is some polynomial and $\mu \in \q$. By Lemma 2.7,
there exists a vector bundle $V$ on $S^{2k-1}\times S^{1}$
satisfying the conditions of Proposition 3.2, such that
$a_{k}(\p(V))\not= 0$. Choose a map $\psi :S^{2k-1} \ra \be $
such that $(\psi \times 1)^{*}\e =V$. We see then that $\psi
^{*}\omega _{k}=a_{k}(\p(V))/[S^{1}]$; moreover, for $1\leq
i\leq k-1$, $\psi ^{*}\omega _{i}\in H^{2i-1}(S^{2k-1})=0$.
Therefore, pulling back equation (1) by $\psi $, we get
$a_{k}(\p(V))/[S^{1}]=\mu \psi ^{*}\theta _{k}$. Since
$a_{k}(\p(V))\not= 0$, we conclude that $\mu \not= 0$, and
dividing by $\mu $, we express $\theta _{k}$ as a polynomial
in $\omega _{1},\ldots ,\omega _{k}$.
{}~~$\Box $

\begin{propn}
If $E$ is a complex vector bundle of rank $n$ on $S^{2}$, then
the cohomology algebra of $\be (S^{2},E)$ is generated by
$a_{i}(\p(\e))/[S^{2}]$, $2\leq i \leq n$.
\end{propn}

{\em Proof}.~~By definition, $\be (S^{2},E)$ is a connected
component of $\Omega ^{2}BU(n)$, which is homotopically
equivalent to $\Omega U(n)$. Therefore (see \cite{PS}, p.68),
$H^{*}(\be )$ is generated by $n-1$ elements $\theta _{i} \in
H^{2i}(\be )$, $1\leq i\leq n-1$. If $\omega
_{i}=a_{i+1}(\p(\e))/[S^{2}]$, $1\leq i \leq n-1$, then using
Lemma 2.7 and Proposition 3.2, we see, as in the proof of
Proposition 3.3, that each $\theta _{k}$ is a polynomial in
$\omega _{1},\ldots ,\omega _{k}$.
{}~~$\Box $

We now use these results to obtain generators for $\be (X,E)$
when $X$ is a 2-manifold.

\begin{propn}
Let $X$ be a pointed, compact, connected and oriented surface,
and let $E$ be a complex vector bundle of rank $n$ on $X$.
Then the cohomology algebra of $\be (X,E)$ is generated by the
images of
$$\sigma (c_{1}(\e)):H_{1}(X)\lra H^{1}(\be )~~~~{\it and}$$

$$\sigma (a_{i}(\p(\e))):H_{r}(X)\ra H^{2i-r}(\be )
{}~~~(2\leq i\leq n,1\leq r\leq 2),$$
where $\sigma $ denotes, as usual, the slant product.
\end{propn}

{\em Proof}.~~Write $X$ as a cofibration
$B\stackrel{i}{\hookrightarrow}X \stackrel{\pi}{\ra} S^{2}$,
where $B$ is a wedge of $2g$ circles, and assume, without loss of
generality, that the centre of the wedge is the base point $x_{0}$ of
$X$. Denote $G=i^{*}E$. Let $z_{0}=\pi (x_{0})$, and let $F$ be a
vector bundle over $S^{2}$ such that $\pi ^{*}F\cong E$. Since the
mapping functor transforms cofibrations into fibrations, we get a
fibration
$$\be (S^{2},F) \stackrel{\pi ^{\#}}{\hookrightarrow}\be (X,E)
\stackrel{i^{\#}}{\ra}\be (B,G).$$
Proposition 3.1 implies that the pull-back of
$a_{i}(\p\e(X,E))/[X]$ by $\pi ^{\#}$ equals
$a_{i}(\p\e(S^{2},F))/\pi _{*}[X]$. Since $H_{2}(\pi )$ is an
isomorphism, Proposition 3.4 now tells us that the Leray-Hirsch
theorem applies. Thus the cohomology algebra of $\be (X,E)$ is
generated by the $a_{i}(\p\e(X,E))/[X]$ together with the image of
$H^{*}(i^{\#})$. Let $B=\vee_{\alpha =1}^{2g}S_{\alpha }$,
where  each $S_{\alpha }$ is a
circle; then $\gamma _{\alpha }=[S_{\alpha }]$ form a basis of
$H_{1}(B)$. Since $\be (B,G)=\prod _{\alpha =1}^{2g}\be (S_{\alpha
},G)$, and since $H_{1}(i)$ is an isomorphism, Propositions 3.1 and
3.3, applied as above, lead us to the finish.
{}~~$\Box $

\begin{thm}
Let $X$ be a 2-manifold as in Proposition 3.5. Then, the cohomology
algebra of $\be\g$ is generated by the images of
$$\sigma (c_{1}(\e)):H_{r}(X)\lra H^{2-r}(\b\g)~~(0\leq r\leq 1)~~{\it
and}$$

$$\sigma (a_{i}(\p(\e))): H_{r}(X) \lra H^{2i-r}(\b\g)~~(0\leq r\leq
2,2\leq i \leq n).$$
\end{thm}

{\em Proof}.~~ Let $x_{0}$ be the base point of $X$, and consider the
fibration $\be \hookrightarrow B\g \stackrel{\varepsilon
_{x_{0}}}{\ra} BU(n)$,
where $\varepsilon _{x_{0}}(f)=f(x_{0})$. Since the images of
$H_{1}(X)$
and $H_{2}(X)$ under the various slant products restrict, by
Proposition 3.5, to generators of $H^{*}(\be )$, the Leray-Hirsch
theorem applies. By Lemma 2.6, $c_{1}(EU(n))$ and $a_{i}(\p EU(n))$,
$2\leq i\leq n$, generate $H^{*}(BU(n))$. Since $\varepsilon
_{x_{0}}^{*}d(EU(n))=d(\e)/[x_{0}]$ for $d=c_{1}(.)$ or
$d=a_{i}(\p(.))$, the result follows.
{}~~$\Box $

\begin{rem}
\begin{em}
In view of Lemma 2.6, Theorem 2.6 implies the assertion concerning
rational cohomology in Proposition 2.20 of \cite{AB}. Actually Lemma
2.6 may give one the impression that the above theorem can be deduced
from Proposition 2.20 of \cite{AB}, but this impression is hard to
substantiate; the difficulty is due to the fact that the slant
product does not behave well with the cup product.
\end{em}
\end{rem}

We now apply the above results in the context of parabolic bundles
over a curve. The standard reference for parabolic bundles is Mehta
and Seshadri \cite{MS}, and we refer to Nitsure \cite{Nit} for the
gauge theoretic aspects of parabolic bundles.

Let $X$ be a compact, connected and oriented surface, fix a positive
integer $n$, and let $\Delta $ be a parabolic datum of rank $n$ on
$X$. Thus, $\Delta $ consists of:

\begin{itemize}
\item a finite subset $J$ of $X$; and \\
\item for each $x\in J$, a sequence $(n_{x,1},\ldots ,n_{x,k_{x}})$
of positive integers such that $\sum _{i=1}^{k_{x}} n_{x,i}=n$,
and a sequence $0\leq \alpha _{x,1}<\ldots <\alpha _{x,k_{x}}<1$ of
real numbers.
\end{itemize}

\noindent Fix a quasi-parabolic vector bundle of rank $n$
and type $\Delta $ on
$X$, and let $\gp $ denote the subgroup of the gauge group $\g $ of
$E$, consisting of parabolic gauge transformations.

Let $\e $ denote the universal bundle on $B\g \times X$, and
for each $x\in J$, let $\f _{x}$ denote the bundle of flags of
type $\Delta $ in $\jxs \e$, where $\jay :B\g \ra B\g \times
X$ is the map $f\mapsto (f,x)$. Define $\phi :\f \ra B\g $ to
be the fibre product of $\f _{x}~~(x\in J)$ over $B\g$. With
these preparations out of the way, we can identify $B\gp $.

\begin{lemma}
The space $\f$ is a classifying space for $\gp$.
\end{lemma}

{\em Proof}.~~General considerations give $B\gp $ as $E\g
/\gp$; further, by \cite{AB}, $E\g$ is the space of all maps
$\tilde{f}:E\ra EU(n)$ which carry each fibre of $E$
isomorphically to some fibre of $EU(n).$ Note that such an
$\tilde{f}$ defines an element $f\in B\g$ such that
$\tilde{f}$ is an isomorphism of $E$ with $f^{*}EU(n)$. On the
other hand, the fibre of $\f$ at $f\in B\g$ is the product of
certain flag manifolds of $EU(n)_{f(x)}~(x\in J)$. Define
$\alpha :E\g \ra \f$ by $\alpha (\tilde{f})=(\tilde{f}
(F^{i}E_{x}) \subset EU(n)_{f(x)})$, where $f \in B\g$ is
induced by $\tilde{f}$, and $F^{i}E_{x}$ are given by the
quasi-parabolic structure of $E$. By the definition of $\gp$,
$\alpha $ factors through a map $\tilde{\alpha }:E\g /\gp \ra
\f$, which is easily seen to be a homeomorphism.~~$\Box$

So, denote $\f$ by $B\gp$. The pull-back $\ep$ of $\e$ by
$\phi \times 1 :B\gp \times X \ra B\g \times X$ is a family of
quasi-parabolic bundles, i.e., for each $x\in J$, there is a
decreasing flag $\exi{1}\supset \exi{2}\supset \ldots $ of
type $\Delta $ in $\jxs \ep$, where, as usual, $\jay
(t)=(t,x)$ for $t\in B\gp$. We call $\ep $ the universal
bundle on $B\gp \times X$.

\begin{thm}
With notation as above, the cohomology algebra of $B\gp $ is
generated by $c_{j}(\homo (\exi{i},\exi{i-1}))~(x\in J)$ and
the images of
$$\sigma (c_{1}(\ep )):H_{r}(X)\ra H^{2-r}(B\gp)~~(0\leq r\leq
1)~~{\it and}$$
$$\sigma (a_{i}(\p\ep)):H_{r}(X)\ra H^{2i-r}(B\gp )~~(0\leq
r\leq 2,~2\leq i\leq n).$$
\end{thm}

{\em Proof}.~~As already remarked, the fibre of $\phi :B\gp
\ra B\g$ over $f$ is a product of flag manifolds $M_{f}^{x}$
of the vector spaces $\e _{(f,x)}~(x\in J)$. Each of these
flag manifolds carries a tautological flag $F_{f}^{x,i}$ of
vector bundles, and the flag $\exi{i}$ on $B\gp$, in fact,
restricts to the tautological flag on each factor $M_{f}^{x}$
of the fibre. Now, in general, if $M$ is a flag manifold, and
if $F^{1}\supset F^{2}\supset \ldots $ is its tautological
flag of vector bundles, then $c_{j}(\homo (F^{i},F^{i-1}))$
generate the cohomology algebra of $M$. In our context, this
fact implies that the Leray-Hirsch theorem holds for the
fibration $\phi $, and the result follows from Theorem 3.6.
{}~~$\Box $

\section{Proofs}

This section brings together the results of the previous sections
to prove Theorems 1.4 and 1.5. The notation is the same as
before.

Let $X$ be a compact Riemann surface, and let $n$ and $d$ be
integers with $n$ positive, and let $\Delta $ be a parabolic
datum of rank $n$ on $X$. {\em Suppose that Assumptions 1.1 and
1.2 are satisfied.}

Let $E$ be a $C^{\infty}$ quasi-parabolic bundle on $X$ of rank
$n$, degree $d$ and parabolic type $\Delta $. Let ${\cal A}$ be the
space of holomorphic structures in $E$, and $\asp $ the open
subset of ${\cal A} $ consisting of holomorphic structures which
are  parabolic stable with respect to the datum $\Delta $. Let
$\g$ denote the gauge group of $E$, and denote $\gb = \g /{\bf
C}^{*}$ and $\gbp = \gp /{\bf C}^{*}$,
where ${\bf C}^{*}$ is the constant scalar subgroup
of $\g$. There is a natural action of $\gb $ on ${\cal A}$, which
induces a free action of $\gbp$ on $\asp$, and there is a
canonical homeomorphism of $\asp /\gbp$ with $\ux$, hence we
will identify them with each other from now on.

\begin{nota}
\begin{em}
If $G$ is a topological group, and $T$ is a $G$-space, then
$T(G)$ denotes the homotopy quotient $EG\times _{G}T$. (We
write $T(G)$ instead of the standard notation $T_{G}$ for
reasons of convenience.)
\end{em}
\end{nota}

\begin{rem}
\begin{em}
Note that for any $G$-space $T$, there are two canonical maps
$T(G) \ra BG$ and $T(G)\ra T/G$. The first map is a fibre
bundle over $BG$ with fibre $T$, and is a homotopy
equivalence if $T$ is contractible. The second map is a
homotopy equivalence of $T(G)$ with $T/G$ if $T$ is a free
$G$-space.
\end{em}
\end{rem}

Consider the diagram

\begin{center}
\def\normalbaselines{\baselineskip17pt \lineskip3pt
\lineskiplimit3pt }
\def\mapright#1{\smash{
    \mathop{\longrightarrow}\limits^{#1}}}
\def\mapdown#1{\Big\downarrow
 \rlap{$\vcenter{\hbox{$\scriptstyle#1$}}$}}

$$\matrix{1&\mapright{}&{\bf
C}^{*}&\mapright{}&{\gp}&\mapright{f}&{\gbp}&\mapright{}&1\cr
&&\mapdown{}&& \mapdown{}&&\mapdown{}&&\cr 1&\mapright{}&{\bf
C}^{*}&\mapright{}& {\g}&\mapright{f'}&\gb &\mapright{}& 1\cr}
$$
\end{center}

\noindent where the vertical maps are the canonical inclusions,
and $f$ and $f'$ denote the canonical projections. This induces a
diagram

\begin{center}
\setlength{\unitlength}{0.5pt}
\begin{picture}(288,210)(0,-20)

\put(20,2){\makebox(0,0){$B\g $}}
\put(256,2){\makebox(0,0){$B\gb $}}
\put(256,162){\makebox(0,0){$B\gbp $}}
\put(20,162){\makebox(0,0){$B\gp $}}

\put(74,2){\vector(1,0){146}}
\put(256,144){\vector(0,-1){124}}
\put(74,162){\vector(1,0){146}}
\put(20,144){\vector(0,-1){124}}

\put(138,-15){\makebox(0,0){$\pi '$}}
\put(138,179){\makebox(0,0){$\pi $}}

\end{picture}
\end{center}
\noindent of fibrations; the fibres of $\pi $ and $\pi '$ are
homeomorphic to $BU(1)$.

Since ${\cal A} $ is contractible, by Remark 4.2, the
natural map\\
$\phi : \ay (\gp)\ra B\gp$ is a homotopy equivalence. Let $\ep $
denote the universal bundle on $B\gp \times X$, and let $V$
denote the bundle on $\asp(\gp) \times X$ obtained by pulling
back $\ep$ via the composition $$\asp (\gp) \times X
\stackrel{\lambda \times 1} {\hookrightarrow}\ay(\gp)\times X
\stackrel{\phi \times 1}{\ra}B\gp \times X,$$
where $\lambda :\asp (\gp)\hookrightarrow \ay(\gp)$ denotes the
inclusion map.

\begin{rem}
\begin{em}
The $\gp$-equivariant perfectness of  a certain
stratification (see Nitsure \cite{Nit}) implies that the
inclusion $\lambda $ induces a surjection in rational
cohomology. Thus, by Theorem 3.9, the Chern classes
$c_{j}(\homo (V^{x,i}, V^{x,i-1}))$ and the slant products
$c_{1}(V)/z~(z\in H_{1}(X))$, $c_{1}(V)/[x_{0}]$ ($x_{0}$ a
fixed base point in $X$) and $a_{i}(\p(V))/[y]~~(y\in
H_{r}(X), ~0\leq r\leq 2,~2\leq i \leq n)$ generate the
algebra $H^{\ast }(\asp (\gp ))$.
\end{em}
\end{rem}

\noindent {\bf Proof of Theorem 1.4:} ~~Let notation be
as above, and
as in Theorem 1.4. Since the action of $\gbp$ on $\asp $ is
free, the canonical map $\psi :\asp (\gbp)\ra \ux $ is a
homotopy equivalence. Let $V'$ denote the bundle on
$\asp(\gp)\times X$ obtained by pulling back the universal
bundle\\ $U\ra \ux \times X$ by the composition
$$ \asp (\gp) \times X \stackrel{\pi \times 1}{\lra} \asp
(\gbp)\times X \stackrel{\psi \times 1}{\lra}\ux \times X$$
where $\pi $ is induced by $\pi :B\gp \ra B\gbp$ above. Recall
now that we have constructed above another family $V$ on $\asp
(\gp )\times X$ using $B\gp$. Now $V$ and $V'$ are families of
parabolic stable bundles parametrized by $\asp (\gp )$ such
that for each $t \in \asp (\gp)$, $V_{t}\cong V'_{t}$.
Therefore, there exists a line bundle $\xi $ on $\asp (\gp)$
such that $V'\cong V\otimes p^{*}\xi $, where $p:\asp
(\gp)\times X \ra \asp (\gp)$ is the canonical projection. This
implies that $\p(V)\cong \p(V')$, $\homo (V^{x,i},V^{x,i-1})
\cong \homo ((V')^{x,i},(V')^{x,i-1})$, and $c_{1}(V)/z =
c_{1}(V')/z$ for all $z\in H_{1}(X)$. Thus, if $W=(\psi \times
1)^{*}U$, then $\p(V)=(\pi \times 1)^{*}\p(W)$, $\homo
(V^{x,i},V^{x,i-1})=\pi ^{*}(\homo (W^{x,i}, W^{x,i-1}))$, and
$c_{1}(V)/z = \pi ^{*}(c_{1}(W)/z)$ for all $z\in H_{1}(X)$.
Further, we easily see that under the inclusion $BU(1)
\hookrightarrow \asp (\g)$ as a fibre of $\pi
,~c_{1}(V)/[x_{0}]$ restricts to a generator of
$H^{2}(BU(1))$. Now in general, if $\pi :E \ra B$ is a
fibration with fibre $F$ such that: (a) the algebra $H^{*}(E)$
is generated by certain classes $\alpha ,\beta _{1},\ldots
,\beta _{k}$; (b) the classes $\beta _{i}$ are pull-backs of
certain classes $\theta _{i} \in H^{\ast}(B)$ by $\pi $; and
(c) $H^{\*}(F)$ is a polynomial algebra on $\alpha _{F}$,
where $\alpha _{F}$ denotes the restriction of $\alpha $
to $F$; then, $H^{*}(B)$ is generated by $\theta _{1},\ldots
,\theta _{k}$. This fact applies in our situation because of
the above observations and because of Remark 4.3, and implies
that the Chern classes $c_{j}(\homo (W^{x,i},W^{x,i-1}))$ and
the slant products $c_{1}(W)/z~(z\in H_{1}(X))$ and
$a_{i}(\p(W))/y~(y\in H_{r}(X),~0\leq r\leq 2,~2\leq i\leq n)$
generate $H^{*}(\asp(\gbp))$. Since $\psi :\asp(\gbp)\ra \ux$
is a homotopy equivalence, we are done. ~~$\Box $

\noindent {\bf Proof of Theorem 1.5:}~~Suppose $U$ is a
universal bundle
on $\sux \times X$. Consider the right action of the
$n$-torsion subgroup $\Gamma _{X}(n)$ of the Jacobian $J_{X}$
on $\sux \times X$, defined by $(E,\alpha ).\zeta = (E\otimes
\zeta ,\zeta ^{-1}\otimes \alpha )$, where $E\in \sux ,~\alpha
\in J_{X}$ and $\zeta \in \Gamma _{X}(n)$. Then the map
$$\pi :\sux \times J_{X} \lra \ux , ~~~(E,\alpha) \mapsto E
\otimes \alpha $$
is a principal $\Gamma
_{X}(n)$-bundle, i.e., a Galois covering with Galois group
$\Gamma _{X}(n)$.   On the other hand, the Poincar\'{e}
polynomials of $\ux$ and\\
$\sux \times J_{X}$ are equal (see Nitsure \cite{Nit},
Remark
3.11). Since $\Gamma _{X}(n)$ is a finite group, this means
that the action of $\Gamma _{X}(n)$ on $\sux \times J_{X}$
induces a trivial action on the cohomology of $\sux \times
J_{X}$, or equivalently that the map $\pi $ induces an
isomorphism in rational cohomology.  (In the case of usual
vector bundles, the triviality of the action of $\gn$ on the
rational cohomology of ${\cal SU}_{X}(n,L)$, where $n$ and the
degree of $L$ are coprime, is a theorem of Harder and
Narasimhan \cite{HN}, who proved it using arithmetic
techniques. It was reproved by Atiyah and Bott \cite{AB} using
gauge theory. The methods of Nitsure \cite{Nit} generalize the
approach of Atiyah and Bott \cite{AB} to parabolic bundles.)
Now, if
$$i:\sux \lra \sux \times J_{X}$$
denotes the map $E\mapsto (E,{\cal O}_{X})$, and if $ j :\sux
\ra \ux $ denotes the inclusion, then $j=\pi \circ i$. Since
$\pi ^{*}$ is an isomorphism and $i^{*}$ is surjective, we see
that $$j^{*}:H^{*}(\ux)\lra H^{*}(\sux )$$ is surjective. Now,
let $\tilde{V}$ be an arbitrary universal bundle on $\ux
\times X$, and denote the restriction of $\tilde{V}$ to $\sux
\times X$ by $V$. Then, Theorem 1.4 applied to $\tilde{V}$,
and the surjectivity of $j^{*}$ imply that the Chern classes
$c_{j}(\homo (V^{x,i},V^{x,i-1}))$ and the slant products
$c_{1}(V)/z~(z\in H_{1}(X))$ and $a_{i}(\p(V))/y~(y\in
H_{r}(X),~0\leq r\leq 2,~2\leq i\leq n)$ generate
$H^{*}(\sux)$. But $\sux $ is simply connected, so the classes
$c_{1}(V)/z$ ($z \in H_{1}(X)$) are all zero. Finally,
since $U$ and $V$ are both
universal bundles on \\ $\sux \times X$, they differ by a line
bundle coming from $\sux$, hence\\ $\homo (V^{x,i},V^{x,i-1})
\cong \homo (U^{x,i},U^{x,i-1})$ and $\p(V)\cong
\p(U)$.~~~~~$\Box $

\noindent
{\bf Proof of Corollary 1.6:}~~If $U^{x}=\jxs U~(x\in J)$,
then the exact sequence
$$0\longrightarrow S^{x} \longrightarrow  U^{x} \longrightarrow
 Q^{x} \longrightarrow 0  $$
implies that $U^{x}\cong S^{x} \oplus Q^{x}$ topologically.
Since $S^{x}$ is either zero or a line bundle,
 $S_{x}^{*} \otimes S_{x}$ is either zero or a trivial line
bundle, and
hence \\ $c_{j}(\homo (S^{x},S^{x}))=c_{j}(\homo (S^{x},Q^{x}))$.
Finally, Example 2.5 (1) implies that \\ $a_{2}(\p(U))=c_{2}(\en
U)$.~~~$\Box $

\noindent
{\bf Proof of Proposition 1.7:}~~As in Atiyah and Bott (see
\cite{AB}, Section 9), the crux of the proof consists in
finding a holomorphic $\gp $-line bundle $\xi $ on $\asp $ on
which ${\bf C}^{*} \subset \gp $ acts via the identity
homomorphism ${\bf C}^{*}\ra {\bf C}^{*},~t \mapsto t$. Let $U
=\asp \times E$ and $U^{x,i}=\asp \times F^{i}E_{x}~(x\in J)$,
where $E$ is the fixed $C^{\infty }$ quasi-parabolic bundle
under consideration. If we let $\gp $ act trivially on $X$,
then $U$ and $U^{x,i}$ are naturally $\gp$-vector bundles on
which ${\bf C}^{*}$ acts by the identity homomorphism. Fix a
line bundle ${\cal O}_{X}(1)$ of degree 1 on $X$, and for each
$k\in {\bf Z}$, let $U(k)=U\otimes q^{*}{\cal O}_{X}(k)$, where
$q:\asp \times X \ra X$ denotes the canonical projection.
Denote by Det $U(k)$ the determinant line bundle of $U(k)$ in
the sense of Quillen \cite{Q}. Then Det $U(k)$ is a
holomorphic $\gp$-line bundle over $\asp $ on which ${\bf
C}^{*}$ acts by the homomorphism $t\mapsto t^{N+kn}$, where
$N=d+n(1-g),~g$ being the genus of $X$. If $(n,d)=1$, let $a,b
\in {\bf Z}$ be such that $an+bN =1$, and take $$\xi = ({\rm
Det}~U(1))^{a}\otimes ({\rm Det}~U)^{b-a};$$ then ${\bf C}^{*}$
acts by the identity homomorphism on $\xi $. If $\sum
_{i=j}^{k_{x}} n_{x,j}$ and $n$ are coprime for some $x$ and
$j$, then the rank $m$ of $U^{x,j}$ and $n$ are coprime; let
$a,b \in {\bf Z}$ be such that $am+bn =1$, and take $$\xi =
({\rm det}~U^{x,j})^{a} \otimes ({\rm Det}~U)^{-b}\otimes
({\rm Det}~U(1))^{b};$$ then $\xi $ has the required property.
Lastly, if $\sum _{i=j}^{k_{x}} n_{x,i}$ and $n+d$ are coprime
for some $x$ and $j$, let $a,b \in {\bf Z}$ be such that $am
+b(d+n)=1$, where $m$ is the rank of  $U^{x,i}$; then $$\xi
=({\rm det}~U^{x,j})^{a}\otimes ({\rm Det}~U(1))^{b} \otimes
({\rm det}~U^{x})^{b(g-1)}$$ will do, where $U^{x}=\jxs
U$.~~~~~$\Box $

\noindent {\em Acknowledgement}.~~ This paper is the outcome
of a suggestion of Professor M.S.Narasimhan. We would like to
thank him for his interest and encouragement.


\begin{thebibliography}{99}
\bibitem{AB} Atiyah, M.F., Bott, R.: The Yang-Mills equations
on Riemann surfaces. Philos. Trans. R. Soc. London {\bf
A 308}, 523-615 (1982)
\bibitem{Beau} Beauville, A.: Sur les cohomologie de certains
espaces du modules de fibr\'{e}s vectoriels. To appear in the
Proceedings of the Int. Coll. on Geometry and Analysis,
Bombay, 1992.
\bibitem{BH} Borel, A., Hirzebruch, F.: Characteristic classes and
homogeneous spaces I. Amer. J. Math. {\bf 80},
458-538 (1958)
\bibitem{DK} Donaldson, S.K., Kronheimer, P.B.:  The
geometry of four-manifolds. Oxford: Clarendon Press 1990
\bibitem{Dre} Drezet, J.M.: Cohomologie du groupe de jauge. In:
Verdier, J-L.,  Le Potier, J. (eds.):
Module des fibr\'{e}s stable sur les courbes alg\'{e}briques, Notes
de l'Ecole Normale Sup\'{e}rieure, Printemps, 1983 (Progress in Math.,
vol. 54, pp. 51-80) Boston Basel Stuttgart : Birkhauser 1985
\bibitem{HN} Harder, G., Narasimhan, M.S.: On the cohomology
groups of moduli spaces. Math. Ann. {\bf 212}, 215-248 (1975)
\bibitem{MS} Mehta, V.B., Seshadri, C.S.: Moduli of vector
bundles on curves with parabolic structures.  Math. Ann.
{\bf 248}, 205-239 (1980)
\bibitem{N}  Newstead, P.E.: Characteristic classes of
stable bundles of rank 2 over an algebraic curve.  Trans.
Am. Math Soc. {\bf 169}, 337-345 (1972)
\bibitem{Nit}  Nitsure, N.: Cohomology of the moduli of parabolic
vector bundles. Proc. Indian Acad. Sci. (Math. Sci) {\bf
95}, 61-77 (1986)
\bibitem{PS}  Pressley, A., Segal, G.:  Loop Groups.
Oxford: Clarendon Press 1986
\bibitem{Q}  Quillen, D.: Determinants of Cauchy-Riemann
operators over a Riemann surface. Funct. Anal. Appl.
{\bf 19}, 31-34 (1985)
\end{thebibliography}
\end{document}